\PassOptionsToPackage{frozencache}{minted}
\documentclass[11pt, a4paper]{amsart}
\usepackage[margin=1.35in]{geometry}

\usepackage[style=alphabetic,
   giveninits=true,
   sorting=nty,
   backend=biber,
   maxbibnames=3,
   isbn=false,
   backref=false]{biblatex}
\DeclareSourcemap{
  \maps[datatype=bibtex]{
    \map{
      \step[fieldsource=doi, final]
      \step[fieldset=url, null]
      \step[fieldset=eprint, null]
      \step[fieldset=eprinttype, null]
      \step[fieldset=archiveprefix, null]
      \step[fieldset=primaryclass, null]
    }
    \map{
      \step[fieldset=language, null]
    }
  }
}
\usepackage{booktabs}
\usepackage[dvipsnames]{xcolor}
\usepackage{longtable}
\usepackage{tikz}
\usepackage{subcaption}
\usepackage{siunitx}
\usepackage{thmtools}
\usepackage[colorlinks=true, pdfstartview=FitV, linkcolor=Blue, citecolor=OliveGreen, urlcolor=BrickRed, breaklinks=false]{hyperref}

\makeatletter
\@ifundefined{newcounteralias}{}{%
  \renewcommand\thmt@autorefsetup{%
    \@xa\def\csname\thmt@envname autorefname\@xa\endcsname
      \@xa{\thmt@thmname}%
  }%
}
\makeatother
\newcommand\scheme[1]{\langle #1 \rangle}

\newcommand\RP{\mathbb{RP}^2}
\newcommand\FF{\mathbb{F}_2}
\newcommand\RR{\mathbb{R}}
\newcommand\dDelta{d\Delta_2}
\newcommand\tri{\mathcal{T}}
\newcommand\cL{\mathcal{L}}
\newcommand\coll{\mathfrak{S}}
\newcommand\splits{\mathcal{H}}
\newcommand\curve[1]{\mathcal{C}(#1)}
\newcommand\signof[1]{\sigma_{#1}}
\newcommand\schemes[1]{R(#1)}
\newcommand\libisotopyurl{\url{https://github.com/polymake/libisotopy}}
\newcommand\splitenumurl{\url{https://github.com/dmg-lab/splitenum}}

\theoremstyle{plain}
\newtheorem{theorem}{Theorem}[section]
\newtheorem{proposition}[theorem]{Proposition}
\newtheorem{lemma}[theorem]{Lemma}
\newtheorem{corollary}[theorem]{Corollary}

\theoremstyle{definition}
\newtheorem{definition}[theorem]{Definition}
\newtheorem{example}[theorem]{Example}
\newtheorem{algorithm}[theorem]{Algorithm}
\theoremstyle{remark}

\title[Complete enumeration of maximal T-curves in degrees $8$ and $9$]{Complete enumeration of maximal T-curves \\ in degrees eight and nine}

\author{Zoe Geiselmann}
\author{Michael Joswig}
\author{Lars Kastner}
\author{Marcel Wack}
\address{Technische Universit\"at Berlin, Institut f\"ur Mathematik, MA 6-2 (IMoS), Fasanenstr.\ 89, 10623 Berlin, Germany}
\email{\{geiselmann,joswig,kastner,wack\}@math.tu-berlin.de}
\address[M.~Joswig]{Max Planck Institute for Mathematics in the Sciences, Leipzig, Germany}

\date{\today}

\keywords{real algebraic curves, combinatorial patchworking, T-curves, Hilbert's sixteenth problem, exhaustive enumeration}

\begin{document}

\begin{abstract}
  We use a highly parallelized backtracking algorithm implemented in CUDA to enumerate every compatible collection of Harnack splits in degrees eight and nine.
  By a theorem of Haas (1997) this yields every maximal T-curve.
  In this way, we find that exactly 38 real schemes of real plane algebraic curves of degree eight with 22 real components and 225 of degree nine with 29 real components are realized by combinatorial patchworking.
\end{abstract}

\maketitle

\section{Introduction}\label{sec:introduction}

The topology of real plane algebraic curves has been studied since Harnack~\cite{1876:Harnack:VieltheiligkeitEbenen}, who showed that a real curve of degree $d$ has at most $\frac{1}{2}(d-1)(d-2)+1$ connected components, and that this bound is attained in every degree. Hilbert~\cite{1891:Hilbert:ReellenZuge} took up this question, and it became the main motivation for the first part of Hilbert's sixteenth problem~\cite{1900:Hilbert:Mathematical}, in which he asked for the classification of real plane algebraic curves up to ambient isotopy of $\RP$. This classification is complete for degree at most seven~\cite{1984:Viro:GluingAlgebraic}, nearly complete for degree eight~\cite{2002:Orevkov:CurvesIsotopy} and wide open from degree nine on.

Many of the constructions behind these classifications were systematized by Viro's patchworking~\cite{1984:Viro:GluingAlgebraic}.
Combinatorial patchworking produces a curve from a regular unimodular triangulation of the lattice points of the triangle $\dDelta$ and a sign at each lattice point.
The resulting curves are called T-curves, and their real schemes can be read off directly from these data.
Itenberg~\cite[Section~7]{1995:Itenberg:RagsdaleCurves} showed that not every real scheme arises in this way.
For sufficiently large degree, there are nonempty real schemes that no T-curve of the same degree realizes, and it was shown in~\cite{2026:GeiselmannEtAl:PatchworkedCurves} that this already happens for $M$-curves of degree eight.
Itenberg and Viro~\cite[p.~22]{1996:ItenbergViro:PatchworkingRagsdale} asked which real schemes of degree eight are realized by T-curves.
Prior to this work, by~\cite{2026:GeiselmannEtAl:PatchworkedCurves}, $38$ of the $89$ real schemes of $M$-curves of degree eight were known to be realized, $12$ were excluded and $39$ remained undecided.

For maximal T-curves there is a structural description discovered by Haas~\cite{1997:Haas:AlgebraicCombinatorial} and revisited in~\cite{2017:BertrandEtAl:HaasRevisited}.
Every maximal T-curve of degree $d$ arises from augmenting the Harnack sign function on a regular unimodular triangulation by a compatible collection of Harnack splits of $\dDelta$, and, conversely, every compatible collection whose segments lie in some regular unimodular triangulation yields a maximal T-curve.
The real schemes of maximal T-curves of a given degree are therefore the images of the regular compatible collections under a map that is cheap to evaluate.
We organize the compatible collections into a search tree and prune the walk through it by two simple reductions: exploiting the symmetry of the triangle and excluding splits that never affect the real scheme.
We then implement a GPU version of the isotopy algorithm of \cite{2026:GeiselmannEtAl:IsotopyComputation} and let the search run on a GPU cluster. The computation took about $23$ minutes for degree~$8$ and about $4$ hours for degree~$9$, with the following results.

\begin{restatable}{theorem}{thmdegreeeight}\label{thm:degree8}
  Exactly $38$ real schemes of $M$-curves of degree $8$ are realized by T-curves.
  In particular, the $38$ real schemes listed in~\cite{2026:GeiselmannEtAl:PatchworkedCurves} are the only ones realized by T-curves.
\end{restatable}

\begin{restatable}{theorem}{thmdegreenine}\label{thm:degree9}
Exactly $225$ real schemes of $M$-curves of degree $9$ are realized by T-curves.
\end{restatable}

The paper is organized as follows. Section~\ref{sec:preliminaries} recalls the necessary background on combinatorial patchworking. Section~\ref{sec:method} describes the backtracking algorithm, and Section~\ref{sec:reductions} provides the theoretical basis for the optimizations. Section~\ref{sec:results} presents the results and details of the computation.

\section{Preliminaries}\label{sec:preliminaries}

We describe combinatorial patchworking following the conventions of~\cite{2026:GeiselmannEtAl:PatchworkedCurves}.
A \emph{real plane algebraic curve} of degree $d$ is the zero set $V_\RR(f) = \{f=0\} \subset \RP$ of a homogeneous polynomial $f \in \RR[x,y,z]$ of degree $d$; we additionally assume that $V_\RR(f)$ is smooth. The connected components of $V_\RR(f)$ are called \emph{loops}, since each of them is the image of an embedding of the circle $\mathbb{S}^1$ into $\RP$. A loop that separates $\RP$ into a disk and a Möbius strip is called an \emph{oval}, and a nonseparating loop is called a \emph{pseudo-line}. By Harnack's theorem~\cite{1876:Harnack:VieltheiligkeitEbenen}, the number of loops is at most $M = \tfrac{1}{2}(d-1)(d-2)+1$, and a curve attaining this bound is called an \emph{$M$-curve} or \emph{maximal}.
The isotopy type of the pair $(\RP, V_\RR(f))$ is called the \emph{real scheme} of the curve; it is determined by the nesting structure of the ovals (and, for odd~$d$, the pseudo-line).

Some real plane algebraic curves can be constructed using Viro's combinatorial patchworking~\cite{1984:Viro:GluingAlgebraic, 2026:GeiselmannEtAl:PatchworkedCurves}.
For a degree $d\geq 1$, let $\dDelta$ be the triangle with vertices $(0,0)$, $(d,0)$, $(0,d)$, and let $A=\dDelta\cap\mathbb{Z}^2$.
Take a unimodular triangulation $\tri$ of $A$ and a sign function $\sigma\colon A\to\FF$.
Reflecting $\tri$ in the coordinate axes yields four symmetric copies; their union, with antipodal boundary points identified, is a triangulation $\tri^\diamond$ of a model of $\RP$.
Extend $\sigma$ to the vertices of $\tri^\diamond$ by
\begin{gather*}
\sigma(-x,y)=\sigma(x,y)+x,\qquad \sigma(x,-y)=\sigma(x,y)+y,\\
\sigma(-x,-y)=\sigma(x,y)+x+y
\end{gather*}
for $(x,y)\in A$, where all values are taken modulo~$2$.
The \emph{patchworked curve} $\curve{\tri,\sigma}$ is the piecewise linear curve obtained by joining, within each triangle of $\tri^\diamond$, the midpoints of those edges whose endpoints have different signs under~$\sigma$.
This makes sense for every unimodular triangulation~$\tri$.
If $\tri$ is regular, then by Viro's combinatorial patchworking theorem~\cite{1984:Viro:GluingAlgebraic}, $\curve{\tri,\sigma}$ is ambient isotopic to a real plane algebraic curve of degree~$d$, called a \emph{T-curve}, which is constructed from the same data $\tri$ and~$\sigma$.
For nonregular $\tri$, the curve $\curve{\tri,\sigma}$ need not be isotopic to an algebraic curve of degree~$d$.

Note that if $\sigma'=\sigma+\ell$ for an affine function $\ell(x,y)=(\epsilon+ix+jy)\bmod 2$ with $\epsilon,i,j\in\FF$, then $\sigma$ and $\sigma'$ are equivalent in the sense of \cite[Section 2.5]{2026:GeiselmannEtAl:PatchworkedCurves}, and the curves $\curve{\tri,\sigma}$ and $\curve{\tri,\sigma'}$ are ambient isotopic.
A distinguished sign function is the \emph{Harnack sign function} given by $\eta(x,y)=(x+1)(y+1)\bmod 2$. The Harnack sign function and the sign functions equivalent to it have the following property. The \emph{parity} of $u=(x,y)\in A$ is $(x,y)\bmod 2\in\FF^2$. In each of the four symmetric copies of $A$, there is a parity $\alpha\in\FF^2$ such that all points of parity~$\alpha$ have one sign under~$\eta$ and all other points have the opposite sign. Furthermore, the isotopy type of a curve $\curve{\tri, \eta}$ does not depend on the choice of the triangulation~$\tri$~\cite{1995:Itenberg:RagsdaleCurves}.

\section{Method}\label{sec:method}

Our method rests on a theorem of Haas, which describes every maximal T-curve of degree $d$ through a compatible collection of Harnack splits.
Both notions come from \cite{2026:GeiselmannEtAl:PatchworkedCurves} and are recalled in Definitions \ref{def:split} and \ref{def:compatible} below.

\begin{theorem}[Haas]\label{thm:haas}
Let $\coll$ be a compatible collection of Harnack splits and let $\tri$ be a regular unimodular triangulation of $A$ containing all segments of the splits in $\coll$.
Then $\curve{\tri,\signof{\coll}}$ is a maximal curve, and its real scheme does not depend on the choice of $\tri$.
Every maximal T-curve of degree $d$ is ambient isotopic to a curve of this form.
\end{theorem}

This is \cite[Theorem 7.3.0.10]{1997:Haas:AlgebraicCombinatorial}, in the form of \cite[Theorem 4.6]{2017:BertrandEtAl:HaasRevisited}.
As explained before Lemma~\ref{lem:tri_independent} below, the real scheme of $\curve{\tri,\signof{\coll}}$ does not depend on $\tri$ even if nonregular triangulations are admitted.
By Theorem~\ref{thm:haas}, determining the real schemes of maximal T-curves of degree $d$ amounts to running through all compatible collections and checking regularity.
We make this precise in Section~\ref{subsec:splits} and show in Lemma~\ref{lem:pairwise} that compatibility is a condition on pairs of splits.
This opens up the enumeration to a backtracking walk (Section~\ref{subsec:backtracking}) that visits every compatible collection exactly once and computes the real scheme of $\curve{\coll}$ at each node.

\subsection{Harnack splits}\label{subsec:splits}

A \emph{primitive segment} is a segment whose endpoints lie in $A$ and which contains no further point of $A$.
Two primitive segments \emph{cross} if they meet in a point that is not a common endpoint.

\begin{definition}[{\cite[Section 3.2]{2026:GeiselmannEtAl:PatchworkedCurves}}]\label{def:split}
A \emph{Harnack split} $S$ is a path in $\dDelta$ made of one primitive segment (a \emph{simple} split) or of two primitive segments sharing an endpoint, the \emph{apex} (a \emph{double} split), such that $S$ meets the boundary of $\dDelta$ exactly in its two endpoints and the points of $S\cap A$ alternate between two distinct parities $\alpha$ and $\beta$.
The split $S$ cuts $\dDelta$ into two \emph{zones} $Z^+$ and $Z^-$ with $Z^+\cap Z^-=S$, where $Z^+$ is the zone containing fewer vertices of the triangle $\dDelta$, and one of the two is fixed once and for all when both contain the same number.
\end{definition}
The surgical twist of \cite[Section 3.4]{2026:GeiselmannEtAl:PatchworkedCurves} changes a sign function along a Harnack split.
We encode it in a function on $A$, so that a collection of splits acts on the sign function by a sum.
\begin{definition}\label{def:twist_mask}
Let $\alpha=(\alpha_1,\alpha_2)$ and $\beta=(\beta_1,\beta_2)$ be the two parities of the points of $S\cap A$ from Definition \ref{def:split}.
The \emph{twist mask} of $S$ is the function $\tau_S\colon A\to\FF$ with $\tau_S(u)=\ell_S(u)$ for $u\in Z^+$ and $\tau_S(u)=0$ otherwise, where $\ell_S\colon \FF^2\to\FF$, given by
\begin{gather*}
\ell_S(x,y)=(\epsilon+ix+jy)\bmod 2,\\
\epsilon=\alpha_1\beta_2+\alpha_2\beta_1,\quad i=\alpha_2+\beta_2,\quad j=\alpha_1+\beta_1,
\end{gather*}
is the unique nonzero affine function that vanishes on $\alpha$ and on $\beta$.
\end{definition}
Adding $\tau_S$ to a sign function is the surgical twist along $S$ of \cite[Section 3.4]{2026:GeiselmannEtAl:PatchworkedCurves}. Replacing $Z^+$ by $Z^-$ in Definition~\ref{def:twist_mask} gives the \emph{opposite twist mask}, which yields an equivalent sign function.

\begin{definition}[{\cite[Section 3.2]{2026:GeiselmannEtAl:PatchworkedCurves}}]\label{def:compatible}
A collection $\coll$ of Harnack splits is \emph{compatible} if some unimodular triangulation of $A$ contains every segment of every split in $\coll$ as an edge, and no two splits in $\coll$ share a segment.
\end{definition}

\begin{lemma}\label{lem:pairwise}
A collection of Harnack splits is compatible if and only if each of its subsets with at most two elements is compatible.
\end{lemma}

\begin{proof}
A subset of a compatible collection is compatible.
Conversely, assume that all subsets of $\coll$ with at most two elements are compatible, and let $E$ be the set of all segments of all splits in $\coll$.
By assumption, the segments in $E$ are pairwise noncrossing, and no two splits in $\coll$ share a segment.
Since $A$ is finite, $E$ is contained in a set $F$ of pairwise noncrossing primitive segments that is maximal under inclusion.
A maximal set of pairwise noncrossing segments on a finite point set in the plane is the edge set of a triangulation of that point set using all its points \cite{2010:DeLoeraRambauSantos:TriangulationsAlgorithms}.
Each triangle of this triangulation $\tri$ contains no point of $A$ besides its three vertices; hence $\tri$ is unimodular.
\end{proof}

For a compatible collection $\coll$, the sign function of Theorem~\ref{thm:haas} is $\signof{\coll}=\eta+\sum_{S\in\coll}\tau_S$.
The sum is taken in $\FF^A$, so it does not depend on the order of the splits in $\coll$, in accordance with \cite[Lemma 15]{2026:GeiselmannEtAl:PatchworkedCurves}.
The following lemma is part of \cite[Lemma 7.2.0.9]{1997:Haas:AlgebraicCombinatorial} and rests on the fact that any two triangulations of a fixed point set in the plane are connected by a sequence of flips \cite[\S3.4.1]{2010:DeLoeraRambauSantos:TriangulationsAlgorithms}.

\begin{lemma}\label{lem:tri_independent}
Let $\coll$ be a compatible collection of Harnack splits, and let $\tri$ and $\tri'$ be unimodular triangulations of $A$, not necessarily regular, that contain all segments of the splits in $\coll$.
Then $\curve{\tri,\signof{\coll}}$ and $\curve{\tri',\signof{\coll}}$ are ambient isotopic.
\end{lemma}

We write $\curve{\coll}$ for $\curve{\tri,\signof{\coll}}$, where $\tri$ is any unimodular triangulation of $A$ containing the segments of the splits in~$\coll$; by Lemma~\ref{lem:tri_independent}, its real scheme is well defined.
We call $\coll$ \emph{regular} if $\tri$ can be chosen regular.
By Theorem~\ref{thm:haas}, the real schemes of maximal T-curves of degree $d$ are exactly the real schemes of the curves $\curve{\coll}$ for regular compatible collections~$\coll$.

\subsection{Backtracking}\label{subsec:backtracking}

We start with a simple tree whose root is the empty collection and whose children are the collections of a single split.
The tree is then recursively created by adding compatible splits.
Given an order on the splits, we create a tree in which each node represents a unique compatible collection of Harnack splits.
Algorithm~\ref{alg:walk} builds this tree in three steps.

\begin{algorithm}\label{alg:walk}\leavevmode
\renewcommand{\labelenumi}{Step \arabic{enumi}.}
\begin{enumerate}
\item \emph{The list of splits.}
We list all Harnack splits $S_1,\dots,S_N$ of $\dDelta$.
For each $S_k$ we store its twist mask $\tau_k=\tau_{S_k}$.

\item \emph{The compatibility graph.}
For each $k$ we store the set $N(k)$ of all $l\neq k$ such that $\{S_k,S_l\}$ is compatible, as a vector of length $N$.
A single split is compatible by the argument in the proof of Lemma~\ref{lem:pairwise}, since its at most two segments meet only at the apex.
Hence, by Lemma~\ref{lem:pairwise}, the compatible collections are exactly the sets $\{S_{k_1},\dots,\allowbreak S_{k_m}\}$ for which $\{k_1,\dots,k_m\}$ is a clique of the graph $\Gamma$ on $\{1,\dots,N\}$ with neighborhoods $N(k)$, the empty set included.

\item \emph{The walk.}
A node of the search tree is an increasing sequence $k_1<\dots<k_m$ together with its candidate set and its sign vector
\[
P=\{\,l>k_m : l\in N(k_1)\cap\dots\cap N(k_m)\,\},\qquad \sigma=\eta+\tau_{k_1}+\dots+\tau_{k_m}.
\]
The root is the empty sequence with $P=\{1,\dots,N\}$ and $\sigma=\eta$.
First, we complete the segments of the collection to form a unimodular triangulation, $\tri$, by greedily adding all primitive segments in a fixed order that do not cross any of the segments chosen so far.
The result is a maximal set of pairwise noncrossing primitive segments, hence a unimodular triangulation by the proof of Lemma~\ref{lem:pairwise}.
This triangulation need not be regular; by Lemma~\ref{lem:tri_independent}, this does not affect the real scheme.
We compute the real scheme of $\curve{\tri,\sigma}$ with the isotopy algorithm of \cite{2026:GeiselmannEtAl:IsotopyComputation} and add it to the result set $R$.
Then we visit, for each $l\in P$ in increasing order, the child $k_1<\dots<k_m<l$ with candidate set $\{\,l'\in P : l'>l,\ l'\in N(l)\,\}$ and sign vector $\sigma+\tau_l$.
\end{enumerate}
\end{algorithm}

\begin{proposition}\label{prop:walk}
The walk of Algorithm~\ref{alg:walk} visits every compatible collection exactly once, and at the node of $\coll$ the sign vector is $\signof{\coll}$.
\end{proposition}

\begin{proof}
The candidate set of the child $k_1<\dots<k_m<l$ is $\{\,l'>l : l'\in N(k_1)\cap\dots\cap N(k_m)\cap N(l)\,\}$; since $l>k_m$, this agrees with the formula for $P$ at that node.
The sign vector of the child is $\eta+\tau_{k_1}+\dots+\tau_{k_m}+\tau_l$, as claimed.
Every node is a clique of $\Gamma$ by induction on $m$. 
The root is the empty clique, and a child adds an index $l\in P$, which is adjacent to all of $k_1,\dots,k_m$.
Every clique $k_1<\dots<k_m$ is visited, by induction on $m$: its prefix $k_1<\dots<k_{m-1}$ is a clique and hence visited, and $k_m$ lies in the candidate set of the prefix, because $k_m>k_{m-1}$ and $k_m$ is adjacent to all of $k_1,\dots,k_{m-1}$.
By Step 2 of Algorithm~\ref{alg:walk} the cliques of $\Gamma$ are exactly the compatible collections.
\end{proof}

\begin{corollary}\label{cor:walk}
After the walk, $R$ is the set of real schemes of the curves $\curve{\coll}$ over all compatible collections~$\coll$.
In particular, $R$ contains the set of real schemes of maximal T-curves of degree $d$, and a scheme in $R$ is the real scheme of a maximal T-curve if and only if it is the real scheme of $\curve{\coll}$ for some regular compatible collection~$\coll$.
\end{corollary}

\begin{proof}
By Proposition~\ref{prop:walk} and Lemma~\ref{lem:tri_independent}, the walk computes the real scheme of $\curve{\coll}$ for every compatible collection $\coll$, and nothing else.
The remaining claims follow from Theorem~\ref{thm:haas}.
\end{proof}

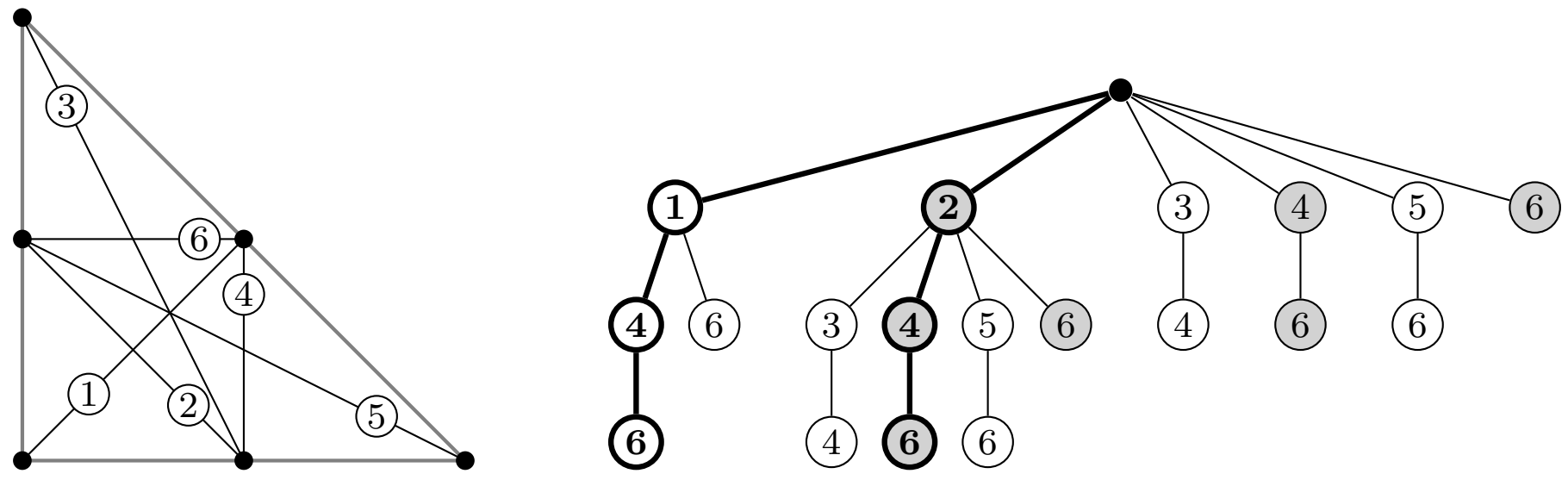
\begin{figure}[t]
\centering
\begin{tikzpicture}[scale=1.7,
  lab/.style={circle,fill=white,draw,inner sep=0.8pt,font=\scriptsize}]
\draw[thick,gray] (0,0)--(2,0)--(0,2)--cycle;
\draw (0,0)--node[lab,pos=0.3]{1}(1,1);
\draw (1,0)--node[lab,pos=0.25]{2}(0,1);
\draw (1,0)--node[lab,pos=0.75]{4}(1,1);
\draw (1,0)--node[lab,pos=0.8]{3}(0,2);
\draw (2,0)--node[lab,pos=0.2]{5}(0,1);
\draw (0,1)--node[lab,pos=0.8]{6}(1,1);
\foreach \p in {(0,0),(1,0),(2,0),(0,1),(1,1),(0,2)} \fill \p circle (1.2pt);
\end{tikzpicture}
\hspace{2em}
\begin{tikzpicture}[x=6mm,y=-9mm,
	t/.style={circle,draw,inner sep=1pt,minimum size=11pt,font=\scriptsize},
	fat/.style={line width=1.2pt,font=\scriptsize\bfseries},
	fatedge/.style={line width=1.2pt},
	gray/.style={fill=gray!35},
	n/a/.style={fat},   n/ad/.style={fat},   n/adf/.style={fat},
	n/b/.style={fat,gray}, n/bd/.style={fat,gray}, n/bdf/.style={fat,gray},
	n/d/.style={gray},  n/f/.style={gray},   n/bf/.style={gray}, n/df/.style={gray},
	e/a/.style={fatedge},  e/b/.style={fatedge},
	e/ad/.style={fatedge}, e/adf/.style={fatedge},
	e/bd/.style={fatedge}, e/bdf/.style={fatedge}]
	\node[circle,fill,inner sep=1.8pt] (r) at (6.2,0) {};
	\foreach \n/\x/\y/\l/\p in {
		a/0.5/1/1/r, b/4/1/2/r, c/7/1/3/r, d/8.5/1/4/r, e/10/1/5/r, f/11.5/1/6/r,
		ad/0/2/4/a, af/1/2/6/a, bc/2.5/2/3/b, bd/3.5/2/4/b, be/4.5/2/5/b, bf/5.5/2/6/b,
		cd/7/2/4/c, df/8.5/2/6/d, ef/10/2/6/e,
		adf/0/3/6/ad, bcd/2.5/3/4/bc, bdf/3.5/3/6/bd, bef/4.5/3/6/be}
	{\node[t,n/\n/.try] (\n) at (\x,\y) {\l}; \draw[e/\n/.try] (\p)--(\n);}
\end{tikzpicture}
\caption{Left: the six Harnack splits of degree $2$. Right: the search tree of Algorithm~\ref{alg:walk}. The root is the empty collection, and a node labeled $l$ adds the split $S_l$ to the collection of its parent. On each level (i.e., distance to the root) nodes of the same color (gray or white) are equivalent under $G$ (see Section~\ref{subsec:symmetry}), where the node with the thick boundary is the canonical node.}
\label{fig:degree2}
\end{figure}

\begin{example}\label{ex:degree2}
For $d=2$ the set $A$ has six points, all on the boundary, so there is no possible apex for a double split and hence every Harnack split is simple.
There are six of them, shown on the left of Figure~\ref{fig:degree2}.
Two of them are compatible exactly when their segments do not cross, which gives the neighborhoods
\begin{gather*}
N(1)=\{4,6\},\quad N(2)=\{3,4,5,6\},\quad N(3)=\{2,4\},\\
N(4)=\{1,2,3,6\},\quad N(5)=\{2,6\},\quad N(6)=\{1,2,4,5\}.
\end{gather*}
The walk of Algorithm~\ref{alg:walk} visits the $20$ nodes on the right of Figure~\ref{fig:degree2}: the root, $6$ collections with one split, $9$ with two and $4$ with three.
For instance, the node $2$ has candidate set $P=N(2)=\{3,4,5,6\}$.
Its child $2<3$ has candidate set $\{\,l'\in P : l'>3,\ l'\in N(3)\,\}=\{4\}$, and its child $2<3<4$ has an empty candidate set.
The four collections with three splits are the maximal ones, and each of them is the set of interior edges of a unimodular triangulation of $A$.
Every node gives a maximal curve by Theorem~\ref{thm:haas}, and the only real scheme of a maximal curve of degree $2$ is one oval, denoted by $\scheme{1}$, so $R=\{\scheme{1}\}$.
\end{example}

The tree has one node per compatible collection, too many to visit in higher degree.

\section{Reductions}\label{sec:reductions}

We shrink the walk of Algorithm~\ref{alg:walk} on the list $\splits$ of all Harnack splits of $\dDelta$ in two ways.
We visit only canonical representatives of the orbits of the symmetry group of the patchwork (Section~\ref{subsec:symmetry}), and we drop splits that never change the real scheme (Section~\ref{subsec:trivial}).
For a list $\cL\subseteq\splits$, let $\schemes{\cL}$ be the set of real schemes of the curves $\curve{\coll}$ over all compatible collections $\coll\subseteq\cL$.
By Corollary~\ref{cor:walk}, $\schemes{\splits}=R$, and the walk of Algorithm~\ref{alg:walk} on $\cL$ computes $\schemes{\cL}$.

\subsection{Symmetry}\label{subsec:symmetry}

The symmetries of combinatorial patchworking have already been discussed in \cite[Section 2.5]{2026:GeiselmannEtAl:PatchworkedCurves}.
Let $G$ be the group of affine unimodular maps of $\RR^2$ that preserve $A$. 
By \cite[Lemma 5]{2026:GeiselmannEtAl:PatchworkedCurves} it is isomorphic to the symmetric group on the three vertices of $\dDelta$.
\begin{proposition}\label{prop:symmetry}
	\begin{enumerate}
		\item For every compatible collection $\coll$ and every $g\in G$, the curves $\curve{\coll}$ and $\curve{g\coll}$ have the same real scheme.
		\item The list $\splits$ is invariant under $G$.
	\end{enumerate}
\end{proposition}
\begin{proof}
	This follows directly from \cite[Lemma 5 and Section 2.5]{2026:GeiselmannEtAl:PatchworkedCurves}.
\end{proof}

We choose one canonical representative to visit for every $G$-orbit of compatible collections. 
Order the list $\splits$ as $S_1,\dots,S_N$, so that $G$ acts on indices, and compare two sets of indices of equal size by their increasing sequences in lexicographic order.
A node of the walk is \emph{canonical} if no element of $G$ maps it to a lexicographically smaller set. 
Observe that if a node $k_1<\dots<k_m$ is canonical, then so is its parent $k_1<\dots<k_{m-1}$, and therefore dismissing noncanonical nodes and their descendants is unproblematic. This proves the following. 

\begin{lemma}\label{lem:symmetry}
	The walk that skips every noncanonical node together with its subtree visits exactly one node in each $G$-orbit of compatible collections and computes the same set of real schemes as the full walk.
\end{lemma}

\begin{example}\label{ex:degree2_cont}
We continue Example~\ref{ex:degree2}.
Under $G$ we have the following equivalences for the index sets of compatible collections. 
\begin{gather*}
	\textbf{\{1\}}\equiv\{3\}\equiv\{5\},\quad \textbf{\{2\}}\equiv\{4\}\equiv\{6\},\\
	\textbf{\{1,4\}}\equiv\{1,6\}\equiv\{2,3\}\equiv\{2,5\}\equiv\{3,4\}\equiv\{5,6\},\\
	\textbf{\{2,4\}}\equiv\{2,6\}\equiv\{4,6\},\quad\textbf{\{1,4,6\}}\equiv \{2,3,4\}\equiv\{2,5,6\},\quad\textbf{\{2,4,6\}}.
\end{gather*}
The canonical representatives are set in bold, and the corresponding canonical nodes are drawn thick in Figure~\ref{fig:degree2}.
In accordance with Lemma~\ref{lem:symmetry}, the thick nodes form a subtree containing the root.
\end{example}

\subsection{Trivial splits}\label{subsec:trivial}

\begin{definition}\label{def:trivial}
A Harnack split $S$ is \emph{trivial} if one of the two zones of $S$ contains no point of $A$ in the interior of $\dDelta$ except possibly the apex of $S$, if $S$ is a double split. This zone is called a \emph{trivial zone}.
\end{definition}

\begin{lemma}\label{lem:ovals_contain_int_points}
  Any oval of a patchworked curve $\curve{\tri,\sigma}$ contains at least one interior lattice point of some symmetric copy of $\dDelta$. 
\end{lemma}

\begin{proof}
  Consider a boundary lattice point $p$ of $\dDelta$; without loss of generality we may assume its second coordinate is zero, so $p=(a,0)$. Since $\tri$ is unimodular, $p$ is adjacent in $\tri$ to at least one point in $A$ of the form $q=(b,1)$. In the patchwork, $p$ is adjacent to both $q$ and $q'=(b,-1)$, where $\sigma(q')=\sigma(q)+1$, so $p$ cannot be separated from both $q$ and $q'$ by $\curve{\tri,\sigma}$. 
\end{proof}

\begin{figure}[t]
	\centering
	\pgfmathsetmacro{\s}{.31*\linewidth/10.5cm}
	\pgfmathsetmacro{\f}{\s/2.646}%
	\definecolor{lightgray}{RGB}{160,160,160}
	\definecolor{darkgray}{RGB}{90,90,90}
	\begin{subfigure}[b]{.32\linewidth}
		\centering
		\begin{tikzpicture}[scale=\s]
			\path (0,-3.25) -- (0,3.25);%
			\begin{scope}
				\draw[lightgray,line width=\f*1.882pt,line cap=round] (0,0) -- (1,0);
				\draw[lightgray,line width=\f*1.882pt,line cap=round] (0,1) -- (1,0);
				\draw[lightgray,line width=\f*1.882pt,line cap=round] (0,0) -- (0,1);
				\draw[lightgray,line width=\f*1.882pt,line cap=round] (0,0) -- (0,-1);
				\draw[lightgray,line width=\f*1.882pt,line cap=round] (1,0) -- (0,-1);
				\draw[lightgray,line width=\f*1.882pt,line cap=round] (0,1) -- (-1,0);
				\draw[lightgray,line width=\f*1.882pt,line cap=round] (0,0) -- (-1,0);
				\draw[lightgray,line width=\f*1.882pt,line cap=round] (0,-1) -- (-1,0);
				\draw[lightgray,line width=\f*1.882pt,line cap=round] (0,1) -- (0,2);
				\draw[lightgray,line width=\f*1.882pt,line cap=round] (0,-1) -- (0,-2);
				\draw[lightgray,line width=\f*1.882pt,line cap=round] (1,0) -- (2,0);
				\draw[lightgray,line width=\f*1.882pt,line cap=round] (0,1) -- (2,0);
				\draw[lightgray,line width=\f*1.882pt,line cap=round] (2,0) -- (0,-1);
				\draw[lightgray,line width=\f*1.882pt,line cap=round] (0,1) -- (-2,0);
				\draw[lightgray,line width=\f*1.882pt,line cap=round] (-1,0) -- (-2,0);
				\draw[lightgray,line width=\f*1.882pt,line cap=round] (0,-1) -- (-2,0);
				\draw[darkgray,line width=\f*4.705pt,line cap=round] (0,1) -- (3,0);
				\draw[darkgray,line width=\f*4.705pt,line cap=round] (3,0) -- (0,-1);
				\draw[darkgray,line width=\f*4.705pt,line cap=round] (0,1) -- (-3,0);
				\draw[darkgray,line width=\f*4.705pt,line cap=round] (0,-1) -- (-3,0);
				\draw[lightgray,line width=\f*1.882pt,line cap=round] (2,0) -- (3,0);
				\draw[lightgray,line width=\f*1.882pt,line cap=round] (-2,0) -- (-3,0);
				\draw[lightgray,line width=\f*1.882pt,line cap=round] (3,0) -- (4,0);
				\draw[lightgray,line width=\f*1.882pt,line cap=round] (-3,0) -- (-4,0);
			\end{scope}
			\begin{scope}
				\filldraw[fill=black, draw=black, line width=\f*2pt] (0,0) circle[radius=\f*15pt];
				\node[anchor=north west, inner sep=1pt, font=\tiny] at (0,0) {$0$};
				\filldraw[fill=black, draw=black, line width=\f*2pt] (0,1) circle[radius=\f*15pt];
				\filldraw[fill=white,draw=black,line width=\f*2pt] (0,2) circle[radius=\f*15pt];
				\filldraw[fill=black, draw=black, line width=\f*2pt] (1,0) circle[radius=\f*15pt];
				\filldraw[fill=black, draw=black, line width=\f*2pt] (2,0) circle[radius=\f*15pt];
				\filldraw[fill=black, draw=black, line width=\f*2pt] (3,0) circle[radius=\f*15pt];
				\filldraw[fill=white,draw=black,line width=\f*2pt] (4,0) circle[radius=\f*15pt];
				\filldraw[fill=white,draw=black,line width=\f*2pt] (0,-1) circle[radius=\f*15pt];
				\filldraw[fill=white,draw=black,line width=\f*2pt] (0,-2) circle[radius=\f*15pt];
				\filldraw[fill=white,draw=black,line width=\f*2pt] (-1,0) circle[radius=\f*15pt];
				\filldraw[fill=black, draw=black, line width=\f*2pt] (-2,0) circle[radius=\f*15pt];
				\filldraw[fill=white,draw=black,line width=\f*2pt] (-3,0) circle[radius=\f*15pt];
				\filldraw[fill=white,draw=black,line width=\f*2pt] (-4,0) circle[radius=\f*15pt];
			\end{scope}
			\begin{scope}
				\draw[black,line width=\f*3.764pt,line cap=round,line join=round] (-0.5,1) -- (-1.5,0.5) -- (-2.5,0) -- (-1,-0.5) -- (-1.5,0) -- (-0.5,0.5) -- (-0.5,0) -- (0,-0.5) -- (0.5,-0.5) -- (1,-0.5) -- (1.5,-0.5) -- (2,-0.5) -- (2.5,-0.5) -- (3,-0.5);
			\end{scope}
		\end{tikzpicture}
		
		\bigskip
		
		\begin{tikzpicture}[scale=\s]
			\path (0,-3.25) -- (0,3.25);%
			\begin{scope}
				\draw[lightgray,line width=\f*1.882pt,line cap=round] (0,0) -- (1,0);
				\draw[lightgray,line width=\f*1.882pt,line cap=round] (0,1) -- (1,0);
				\draw[lightgray,line width=\f*1.882pt,line cap=round] (0,0) -- (0,1);
				\draw[lightgray,line width=\f*1.882pt,line cap=round] (0,0) -- (0,-1);
				\draw[lightgray,line width=\f*1.882pt,line cap=round] (1,0) -- (0,-1);
				\draw[lightgray,line width=\f*1.882pt,line cap=round] (0,1) -- (-1,0);
				\draw[lightgray,line width=\f*1.882pt,line cap=round] (0,0) -- (-1,0);
				\draw[lightgray,line width=\f*1.882pt,line cap=round] (0,-1) -- (-1,0);
				\draw[lightgray,line width=\f*1.882pt,line cap=round] (0,1) -- (0,2);
				\draw[lightgray,line width=\f*1.882pt,line cap=round] (0,-1) -- (0,-2);
				\draw[lightgray,line width=\f*1.882pt,line cap=round] (1,0) -- (2,0);
				\draw[lightgray,line width=\f*1.882pt,line cap=round] (0,1) -- (2,0);
				\draw[lightgray,line width=\f*1.882pt,line cap=round] (2,0) -- (0,-1);
				\draw[lightgray,line width=\f*1.882pt,line cap=round] (0,1) -- (-2,0);
				\draw[lightgray,line width=\f*1.882pt,line cap=round] (-1,0) -- (-2,0);
				\draw[lightgray,line width=\f*1.882pt,line cap=round] (0,-1) -- (-2,0);
				\draw[darkgray,line width=\f*4.705pt,line cap=round] (0,1) -- (3,0);
				\draw[darkgray,line width=\f*4.705pt,line cap=round] (3,0) -- (0,-1);
				\draw[darkgray,line width=\f*4.705pt,line cap=round] (0,1) -- (-3,0);
				\draw[darkgray,line width=\f*4.705pt,line cap=round] (0,-1) -- (-3,0);
				\draw[lightgray,line width=\f*1.882pt,line cap=round] (2,0) -- (3,0);
				\draw[lightgray,line width=\f*1.882pt,line cap=round] (-2,0) -- (-3,0);
				\draw[lightgray,line width=\f*1.882pt,line cap=round] (3,0) -- (4,0);
				\draw[lightgray,line width=\f*1.882pt,line cap=round] (-3,0) -- (-4,0);
			\end{scope}
			\begin{scope}
				\filldraw[fill=white,draw=black,line width=\f*2pt] (0,0) circle[radius=\f*15pt];
				\node[anchor=north west, inner sep=1pt, font=\tiny] at (0,0) {$0$};
				\filldraw[fill=black, draw=black, line width=\f*2pt] (0,1) circle[radius=\f*15pt];
				\filldraw[fill=white,draw=black,line width=\f*2pt] (0,2) circle[radius=\f*15pt];
				\filldraw[fill=black, draw=black, line width=\f*2pt] (1,0) circle[radius=\f*15pt];
				\filldraw[fill=white,draw=black,line width=\f*2pt] (2,0) circle[radius=\f*15pt];
				\filldraw[fill=black, draw=black, line width=\f*2pt] (3,0) circle[radius=\f*15pt];
				\filldraw[fill=white,draw=black,line width=\f*2pt] (4,0) circle[radius=\f*15pt];
				\filldraw[fill=white,draw=black,line width=\f*2pt] (0,-1) circle[radius=\f*15pt];
				\filldraw[fill=white,draw=black,line width=\f*2pt] (0,-2) circle[radius=\f*15pt];
				\filldraw[fill=white,draw=black,line width=\f*2pt] (-1,0) circle[radius=\f*15pt];
				\filldraw[fill=white,draw=black,line width=\f*2pt] (-2,0) circle[radius=\f*15pt];
				\filldraw[fill=white,draw=black,line width=\f*2pt] (-3,0) circle[radius=\f*15pt];
				\filldraw[fill=white,draw=black,line width=\f*2pt] (-4,0) circle[radius=\f*15pt];
			\end{scope}
			\begin{scope}
				\draw[black,line width=\f*3.764pt,line cap=round,line join=round] (-0.5,1) -- (-1.5,0.5) -- (-1,0.5) -- (-0.5,0.5) -- (0,0.5) -- (0.5,0) -- (0.5,-0.5) -- (1.5,0) -- (1,0.5) -- (2.5,0) -- (1.5,-0.5) -- (2,-0.5) -- (2.5,-0.5) -- (3,-0.5);
			\end{scope}
		\end{tikzpicture}
		\caption{$S$ is a simple split.}
		\label{fig:simple}
	\end{subfigure}%
	\hfill
	\begin{subfigure}[b]{.32\linewidth}
		\centering
		\begin{tikzpicture}[scale=\s]
			\path (0,-3.25) -- (0,3.25);%
			\begin{scope}
				\draw[lightgray,line width=\f*1.882pt,line cap=round] (0,0) -- (0,1);
				\draw[lightgray,line width=\f*1.882pt,line cap=round] (0,0) -- (0,-1);
				\draw[lightgray,line width=\f*1.882pt,line cap=round] (0,0) -- (1,0);
				\draw[lightgray,line width=\f*1.882pt,line cap=round] (1,0) -- (3,1);
				\draw[darkgray,line width=\f*4.705pt,line cap=round] (0,0) -- (3,1);
				\draw[darkgray,line width=\f*4.705pt,line cap=round] (0,0) -- (3,-1);
				\draw[lightgray,line width=\f*1.882pt,line cap=round] (1,0) -- (3,-1);
				\draw[darkgray,line width=\f*4.705pt,line cap=round] (0,0) -- (-3,1);
				\draw[lightgray,line width=\f*1.882pt,line cap=round] (-1,0) -- (-3,1);
				\draw[lightgray,line width=\f*1.882pt,line cap=round] (0,0) -- (-1,0);
				\draw[lightgray,line width=\f*1.882pt,line cap=round] (-1,0) -- (-3,-1);
				\draw[darkgray,line width=\f*4.705pt,line cap=round] (0,0) -- (-3,-1);
				\draw[lightgray,line width=\f*1.882pt,line cap=round] (1,0) -- (2,0);
				\draw[lightgray,line width=\f*1.882pt,line cap=round] (2,0) -- (3,1);
				\draw[lightgray,line width=\f*1.882pt,line cap=round] (2,0) -- (3,-1);
				\draw[lightgray,line width=\f*1.882pt,line cap=round] (-2,0) -- (-3,1);
				\draw[lightgray,line width=\f*1.882pt,line cap=round] (-1,0) -- (-2,0);
				\draw[lightgray,line width=\f*1.882pt,line cap=round] (-2,0) -- (-3,-1);
				\draw[lightgray,line width=\f*1.882pt,line cap=round] (2,0) -- (3,0);
				\draw[lightgray,line width=\f*1.882pt,line cap=round] (3,0) -- (3,1);
				\draw[lightgray,line width=\f*1.882pt,line cap=round] (3,0) -- (3,-1);
				\draw[lightgray,line width=\f*1.882pt,line cap=round] (-3,0) -- (-3,1);
				\draw[lightgray,line width=\f*1.882pt,line cap=round] (-2,0) -- (-3,0);
				\draw[lightgray,line width=\f*1.882pt,line cap=round] (-3,0) -- (-3,-1);
				\draw[lightgray,line width=\f*1.882pt,line cap=round] (3,0) -- (4,0);
				\draw[darkgray,line width=\f*4.705pt,line cap=round] (3,1) -- (4,0);
				\draw[darkgray,line width=\f*4.705pt,line cap=round] (4,0) -- (3,-1);
				\draw[darkgray,line width=\f*4.705pt,line cap=round] (-3,1) -- (-4,0);
				\draw[lightgray,line width=\f*1.882pt,line cap=round] (-3,0) -- (-4,0);
				\draw[darkgray,line width=\f*4.705pt,line cap=round] (-3,-1) -- (-4,0);
				\draw[lightgray,line width=\f*1.882pt,line cap=round] (4,0) -- (5,0);
				\draw[lightgray,line width=\f*1.882pt,line cap=round] (-4,0) -- (-5,0);
			\end{scope}
			\begin{scope}
				\filldraw[fill=white,draw=black,line width=\f*2pt] (0,0) circle[radius=\f*15pt];
				\node[anchor=north west, inner sep=1pt, font=\tiny] at (0,-0.15) {$0$};
				\filldraw[fill=black, draw=black, line width=\f*2pt] (0,1) circle[radius=\f*15pt];
				\filldraw[fill=white,draw=black,line width=\f*2pt] (1,0) circle[radius=\f*15pt];
				\filldraw[fill=white,draw=black,line width=\f*2pt] (2,0) circle[radius=\f*15pt];
				\filldraw[fill=white,draw=black,line width=\f*2pt] (3,0) circle[radius=\f*15pt];
				\filldraw[fill=black, draw=black, line width=\f*2pt] (3,1) circle[radius=\f*15pt];
				\filldraw[fill=white,draw=black,line width=\f*2pt] (4,0) circle[radius=\f*15pt];
				\filldraw[fill=black, draw=black, line width=\f*2pt] (5,0) circle[radius=\f*15pt];
				\filldraw[fill=white,draw=black,line width=\f*2pt] (0,-1) circle[radius=\f*15pt];
				\filldraw[fill=white,draw=black,line width=\f*2pt] (3,-1) circle[radius=\f*15pt];
				\filldraw[fill=black, draw=black, line width=\f*2pt] (-1,0) circle[radius=\f*15pt];
				\filldraw[fill=white,draw=black,line width=\f*2pt] (-2,0) circle[radius=\f*15pt];
				\filldraw[fill=black, draw=black, line width=\f*2pt] (-3,0) circle[radius=\f*15pt];
				\filldraw[fill=black, draw=black, line width=\f*2pt] (-3,-1) circle[radius=\f*15pt];
				\filldraw[fill=white,draw=black,line width=\f*2pt] (-3,1) circle[radius=\f*15pt];
				\filldraw[fill=white,draw=black,line width=\f*2pt] (-4,0) circle[radius=\f*15pt];
				\filldraw[fill=white,draw=black,line width=\f*2pt] (-5,0) circle[radius=\f*15pt];
			\end{scope}
			\begin{scope}
				\draw[black,line width=\f*3.764pt,line cap=round,line join=round] (1,0.5) -- (1.5,0.5) -- (2,0.5) -- (2.5,0.5) -- (3,0.5) -- (3.5,0.5) -- (4,0.5);
				\draw[black,line width=\f*3.764pt,line cap=round,line join=round] (-3.5,-1) -- (-3.5,-0.5) -- (-3.5,0) -- (-3,0.5) -- (-2.5,0) -- (-2.5,-0.5) -- (-1.5,0) -- (-2,0.5) -- (-0.5,0) -- (-1.5,-0.5) -- (-2.5,-1);
			\end{scope}
		\end{tikzpicture}
		
		\bigskip
		
		\begin{tikzpicture}[scale=\s]
			\path (0,-3.25) -- (0,3.25);%
			\begin{scope}
				\draw[lightgray,line width=\f*1.882pt,line cap=round] (0,0) -- (0,1);
				\draw[lightgray,line width=\f*1.882pt,line cap=round] (0,0) -- (0,-1);
				\draw[lightgray,line width=\f*1.882pt,line cap=round] (0,0) -- (1,0);
				\draw[lightgray,line width=\f*1.882pt,line cap=round] (1,0) -- (3,1);
				\draw[darkgray,line width=\f*4.705pt,line cap=round] (0,0) -- (3,1);
				\draw[darkgray,line width=\f*4.705pt,line cap=round] (0,0) -- (3,-1);
				\draw[lightgray,line width=\f*1.882pt,line cap=round] (1,0) -- (3,-1);
				\draw[darkgray,line width=\f*4.705pt,line cap=round] (0,0) -- (-3,1);
				\draw[lightgray,line width=\f*1.882pt,line cap=round] (-1,0) -- (-3,1);
				\draw[lightgray,line width=\f*1.882pt,line cap=round] (0,0) -- (-1,0);
				\draw[lightgray,line width=\f*1.882pt,line cap=round] (-1,0) -- (-3,-1);
				\draw[darkgray,line width=\f*4.705pt,line cap=round] (0,0) -- (-3,-1);
				\draw[lightgray,line width=\f*1.882pt,line cap=round] (1,0) -- (2,0);
				\draw[lightgray,line width=\f*1.882pt,line cap=round] (2,0) -- (3,1);
				\draw[lightgray,line width=\f*1.882pt,line cap=round] (2,0) -- (3,-1);
				\draw[lightgray,line width=\f*1.882pt,line cap=round] (-2,0) -- (-3,1);
				\draw[lightgray,line width=\f*1.882pt,line cap=round] (-1,0) -- (-2,0);
				\draw[lightgray,line width=\f*1.882pt,line cap=round] (-2,0) -- (-3,-1);
				\draw[lightgray,line width=\f*1.882pt,line cap=round] (2,0) -- (3,0);
				\draw[lightgray,line width=\f*1.882pt,line cap=round] (3,0) -- (3,1);
				\draw[lightgray,line width=\f*1.882pt,line cap=round] (3,0) -- (3,-1);
				\draw[lightgray,line width=\f*1.882pt,line cap=round] (-3,0) -- (-3,1);
				\draw[lightgray,line width=\f*1.882pt,line cap=round] (-2,0) -- (-3,0);
				\draw[lightgray,line width=\f*1.882pt,line cap=round] (-3,0) -- (-3,-1);
				\draw[lightgray,line width=\f*1.882pt,line cap=round] (3,0) -- (4,0);
				\draw[darkgray,line width=\f*4.705pt,line cap=round] (3,1) -- (4,0);
				\draw[darkgray,line width=\f*4.705pt,line cap=round] (4,0) -- (3,-1);
				\draw[darkgray,line width=\f*4.705pt,line cap=round] (-3,1) -- (-4,0);
				\draw[lightgray,line width=\f*1.882pt,line cap=round] (-3,0) -- (-4,0);
				\draw[darkgray,line width=\f*4.705pt,line cap=round] (-3,-1) -- (-4,0);
				\draw[lightgray,line width=\f*1.882pt,line cap=round] (4,0) -- (5,0);
				\draw[lightgray,line width=\f*1.882pt,line cap=round] (-4,0) -- (-5,0);
			\end{scope}
			\begin{scope}
				\filldraw[fill=white,draw=black,line width=\f*2pt] (0,0) circle[radius=\f*15pt];
				\node[anchor=north west, inner sep=1pt, font=\tiny] at (0,-0.15) {$0$};
				\filldraw[fill=black, draw=black, line width=\f*2pt] (0,1) circle[radius=\f*15pt];
				\filldraw[fill=black, draw=black, line width=\f*2pt] (1,0) circle[radius=\f*15pt];
				\filldraw[fill=white,draw=black,line width=\f*2pt] (2,0) circle[radius=\f*15pt];
				\filldraw[fill=black, draw=black, line width=\f*2pt] (3,0) circle[radius=\f*15pt];
				\filldraw[fill=black, draw=black, line width=\f*2pt] (3,1) circle[radius=\f*15pt];
				\filldraw[fill=white,draw=black,line width=\f*2pt] (4,0) circle[radius=\f*15pt];
				\filldraw[fill=black, draw=black, line width=\f*2pt] (5,0) circle[radius=\f*15pt];
				\filldraw[fill=white,draw=black,line width=\f*2pt] (0,-1) circle[radius=\f*15pt];
				\filldraw[fill=white,draw=black,line width=\f*2pt] (3,-1) circle[radius=\f*15pt];
				\filldraw[fill=white,draw=black,line width=\f*2pt] (-1,0) circle[radius=\f*15pt];
				\filldraw[fill=white,draw=black,line width=\f*2pt] (-2,0) circle[radius=\f*15pt];
				\filldraw[fill=white,draw=black,line width=\f*2pt] (-3,0) circle[radius=\f*15pt];
				\filldraw[fill=black, draw=black, line width=\f*2pt] (-3,-1) circle[radius=\f*15pt];
				\filldraw[fill=white,draw=black,line width=\f*2pt] (-3,1) circle[radius=\f*15pt];
				\filldraw[fill=white,draw=black,line width=\f*2pt] (-4,0) circle[radius=\f*15pt];
				\filldraw[fill=white,draw=black,line width=\f*2pt] (-5,0) circle[radius=\f*15pt];
			\end{scope}
			\begin{scope}
				\draw[black,line width=\f*3.764pt,line cap=round,line join=round] (1,0.5) -- (1.5,0.5) -- (0.5,0) -- (2,-0.5) -- (1.5,0) -- (2.5,0.5) -- (2.5,0) -- (3,-0.5) -- (3.5,0) -- (3.5,0.5) -- (4,0.5);
				\draw[black,line width=\f*3.764pt,line cap=round,line join=round] (-3.5,-1) -- (-3.5,-0.5) -- (-3,-0.5) -- (-2.5,-0.5) -- (-2,-0.5) -- (-1.5,-0.5) -- (-2.5,-1);
			\end{scope}
		\end{tikzpicture}
		\caption{$S$ is a one-edge double split.}
		\label{fig:double1}
	\end{subfigure}%
	\hfill
	\begin{subfigure}[b]{.32\linewidth}
		\centering
		\begin{tikzpicture}[scale=\s]
			\path (0,-3.25) -- (0,3.25);%
			\begin{scope}
				\draw[lightgray,line width=\f*1.882pt,line cap=round] (0,0) -- (1,1);
				\draw[lightgray,line width=\f*1.882pt,line cap=round] (0,1) -- (1,1);
				\draw[lightgray,line width=\f*1.882pt,line cap=round] (0,0) -- (0,1);
				\draw[lightgray,line width=\f*1.882pt,line cap=round] (0,0) -- (0,-1);
				\draw[lightgray,line width=\f*1.882pt,line cap=round] (0,-1) -- (1,-1);
				\draw[lightgray,line width=\f*1.882pt,line cap=round] (0,0) -- (1,-1);
				\draw[lightgray,line width=\f*1.882pt,line cap=round] (0,1) -- (-1,1);
				\draw[lightgray,line width=\f*1.882pt,line cap=round] (0,0) -- (-1,1);
				\draw[lightgray,line width=\f*1.882pt,line cap=round] (0,0) -- (-1,-1);
				\draw[lightgray,line width=\f*1.882pt,line cap=round] (0,-1) -- (-1,-1);
				\draw[lightgray,line width=\f*1.882pt,line cap=round] (0,0) -- (1,0);
				\draw[lightgray,line width=\f*1.882pt,line cap=round] (1,0) -- (1,1);
				\draw[lightgray,line width=\f*1.882pt,line cap=round] (1,0) -- (1,-1);
				\draw[lightgray,line width=\f*1.882pt,line cap=round] (-1,0) -- (-1,1);
				\draw[lightgray,line width=\f*1.882pt,line cap=round] (0,0) -- (-1,0);
				\draw[lightgray,line width=\f*1.882pt,line cap=round] (-1,0) -- (-1,-1);
				\draw[darkgray,line width=\f*4.705pt,line cap=round] (0,2) -- (1,1);
				\draw[lightgray,line width=\f*1.882pt,line cap=round] (0,1) -- (0,2);
				\draw[lightgray,line width=\f*1.882pt,line cap=round] (0,-1) -- (0,-2);
				\draw[darkgray,line width=\f*4.705pt,line cap=round] (0,-2) -- (1,-1);
				\draw[darkgray,line width=\f*4.705pt,line cap=round] (0,2) -- (-1,1);
				\draw[darkgray,line width=\f*4.705pt,line cap=round] (0,-2) -- (-1,-1);
				\draw[lightgray,line width=\f*1.882pt,line cap=round] (0,2) -- (0,3);
				\draw[lightgray,line width=\f*1.882pt,line cap=round] (0,-2) -- (0,-3);
				\draw[lightgray,line width=\f*1.882pt,line cap=round] (1,0) -- (2,0);
				\draw[lightgray,line width=\f*1.882pt,line cap=round] (1,1) -- (2,0);
				\draw[lightgray,line width=\f*1.882pt,line cap=round] (2,0) -- (1,-1);
				\draw[lightgray,line width=\f*1.882pt,line cap=round] (-1,1) -- (-2,0);
				\draw[lightgray,line width=\f*1.882pt,line cap=round] (-1,0) -- (-2,0);
				\draw[lightgray,line width=\f*1.882pt,line cap=round] (-1,-1) -- (-2,0);
				\draw[lightgray,line width=\f*1.882pt,line cap=round] (2,0) -- (3,0);
				\draw[lightgray,line width=\f*1.882pt,line cap=round] (1,1) -- (3,0);
				\draw[lightgray,line width=\f*1.882pt,line cap=round] (3,0) -- (1,-1);
				\draw[lightgray,line width=\f*1.882pt,line cap=round] (-1,1) -- (-3,0);
				\draw[lightgray,line width=\f*1.882pt,line cap=round] (-2,0) -- (-3,0);
				\draw[lightgray,line width=\f*1.882pt,line cap=round] (-1,-1) -- (-3,0);
				\draw[darkgray,line width=\f*4.705pt,line cap=round] (1,1) -- (4,0);
				\draw[darkgray,line width=\f*4.705pt,line cap=round] (4,0) -- (1,-1);
				\draw[darkgray,line width=\f*4.705pt,line cap=round] (-1,1) -- (-4,0);
				\draw[darkgray,line width=\f*4.705pt,line cap=round] (-1,-1) -- (-4,0);
				\draw[lightgray,line width=\f*1.882pt,line cap=round] (3,0) -- (4,0);
				\draw[lightgray,line width=\f*1.882pt,line cap=round] (-3,0) -- (-4,0);
				\draw[lightgray,line width=\f*1.882pt,line cap=round] (4,0) -- (5,0);
				\draw[lightgray,line width=\f*1.882pt,line cap=round] (-4,0) -- (-5,0);
			\end{scope}
			\begin{scope}
				\filldraw[fill=white,draw=black,line width=\f*2pt] (0,0) circle[radius=\f*15pt];
				\node[anchor=north west, inner sep=1pt, font=\tiny] at (0,0) {$0$};
				\filldraw[fill=white,draw=black,line width=\f*2pt] (0,1) circle[radius=\f*15pt];
				\filldraw[fill=white,draw=black,line width=\f*2pt] (0,2) circle[radius=\f*15pt];
				\filldraw[fill=black, draw=black, line width=\f*2pt] (0,3) circle[radius=\f*15pt];
				\filldraw[fill=white,draw=black,line width=\f*2pt] (1,0) circle[radius=\f*15pt];
				\filldraw[fill=black, draw=black, line width=\f*2pt] (1,1) circle[radius=\f*15pt];
				\node[anchor=south west, inner sep=1pt, font=\tiny] at (1,1) {$p$};
				\filldraw[fill=white,draw=black,line width=\f*2pt] (2,0) circle[radius=\f*15pt];
				\filldraw[fill=white,draw=black,line width=\f*2pt] (3,0) circle[radius=\f*15pt];
				\filldraw[fill=white,draw=black,line width=\f*2pt] (4,0) circle[radius=\f*15pt];
				\filldraw[fill=black, draw=black, line width=\f*2pt] (5,0) circle[radius=\f*15pt];
				\filldraw[fill=black, draw=black, line width=\f*2pt] (0,-1) circle[radius=\f*15pt];
				\filldraw[fill=white,draw=black,line width=\f*2pt] (0,-2) circle[radius=\f*15pt];
				\filldraw[fill=white,draw=black,line width=\f*2pt] (0,-3) circle[radius=\f*15pt];
				\filldraw[fill=white,draw=black,line width=\f*2pt] (1,-1) circle[radius=\f*15pt];
				\filldraw[fill=black, draw=black, line width=\f*2pt] (-1,0) circle[radius=\f*15pt];
				\filldraw[fill=black, draw=black, line width=\f*2pt] (-1,-1) circle[radius=\f*15pt];
				\filldraw[fill=white,draw=black,line width=\f*2pt] (-1,1) circle[radius=\f*15pt];
				\filldraw[fill=white,draw=black,line width=\f*2pt] (-2,0) circle[radius=\f*15pt];
				\filldraw[fill=black, draw=black, line width=\f*2pt] (-3,0) circle[radius=\f*15pt];
				\filldraw[fill=white,draw=black,line width=\f*2pt] (-4,0) circle[radius=\f*15pt];
				\filldraw[fill=white,draw=black,line width=\f*2pt] (-5,0) circle[radius=\f*15pt];
			\end{scope}
			\begin{scope}
				\draw[black,line width=\f*3.764pt,line cap=round,line join=round] (0.5,2) -- (0.5,1.5) -- (0.5,1) -- (0.5,0.5) -- (1,0.5) -- (1.5,0.5) -- (2,0.5) -- (2.5,0.5) -- (3,0.5);
				\draw[black,line width=\f*3.764pt,line cap=round,line join=round] (-1.5,-1) -- (-2.5,-0.5) -- (-3.5,0) -- (-2,0.5) -- (-2.5,0) -- (-1.5,-0.5) -- (-1.5,0) -- (-1,0.5) -- (-0.5,0) -- (-0.5,-0.5) -- (0,-0.5) -- (0.5,-1) -- (0,-1.5) -- (-0.5,-1.5) -- (-1,-1.5);
			\end{scope}
		\end{tikzpicture}
		
		\bigskip
		
		\begin{tikzpicture}[scale=\s]
			\path (0,-3.25) -- (0,3.25);%
			\begin{scope}
				\draw[lightgray,line width=\f*1.882pt,line cap=round] (0,0) -- (1,1);
				\draw[lightgray,line width=\f*1.882pt,line cap=round] (0,1) -- (1,1);
				\draw[lightgray,line width=\f*1.882pt,line cap=round] (0,0) -- (0,1);
				\draw[lightgray,line width=\f*1.882pt,line cap=round] (0,0) -- (0,-1);
				\draw[lightgray,line width=\f*1.882pt,line cap=round] (0,-1) -- (1,-1);
				\draw[lightgray,line width=\f*1.882pt,line cap=round] (0,0) -- (1,-1);
				\draw[lightgray,line width=\f*1.882pt,line cap=round] (0,1) -- (-1,1);
				\draw[lightgray,line width=\f*1.882pt,line cap=round] (0,0) -- (-1,1);
				\draw[lightgray,line width=\f*1.882pt,line cap=round] (0,0) -- (-1,-1);
				\draw[lightgray,line width=\f*1.882pt,line cap=round] (0,-1) -- (-1,-1);
				\draw[lightgray,line width=\f*1.882pt,line cap=round] (0,0) -- (1,0);
				\draw[lightgray,line width=\f*1.882pt,line cap=round] (1,0) -- (1,1);
				\draw[lightgray,line width=\f*1.882pt,line cap=round] (1,0) -- (1,-1);
				\draw[lightgray,line width=\f*1.882pt,line cap=round] (-1,0) -- (-1,1);
				\draw[lightgray,line width=\f*1.882pt,line cap=round] (0,0) -- (-1,0);
				\draw[lightgray,line width=\f*1.882pt,line cap=round] (-1,0) -- (-1,-1);
				\draw[darkgray,line width=\f*4.705pt,line cap=round] (0,2) -- (1,1);
				\draw[lightgray,line width=\f*1.882pt,line cap=round] (0,1) -- (0,2);
				\draw[lightgray,line width=\f*1.882pt,line cap=round] (0,-1) -- (0,-2);
				\draw[darkgray,line width=\f*4.705pt,line cap=round] (0,-2) -- (1,-1);
				\draw[darkgray,line width=\f*4.705pt,line cap=round] (0,2) -- (-1,1);
				\draw[darkgray,line width=\f*4.705pt,line cap=round] (0,-2) -- (-1,-1);
				\draw[lightgray,line width=\f*1.882pt,line cap=round] (0,2) -- (0,3);
				\draw[lightgray,line width=\f*1.882pt,line cap=round] (0,-2) -- (0,-3);
				\draw[lightgray,line width=\f*1.882pt,line cap=round] (1,0) -- (2,0);
				\draw[lightgray,line width=\f*1.882pt,line cap=round] (1,1) -- (2,0);
				\draw[lightgray,line width=\f*1.882pt,line cap=round] (2,0) -- (1,-1);
				\draw[lightgray,line width=\f*1.882pt,line cap=round] (-1,1) -- (-2,0);
				\draw[lightgray,line width=\f*1.882pt,line cap=round] (-1,0) -- (-2,0);
				\draw[lightgray,line width=\f*1.882pt,line cap=round] (-1,-1) -- (-2,0);
				\draw[lightgray,line width=\f*1.882pt,line cap=round] (2,0) -- (3,0);
				\draw[lightgray,line width=\f*1.882pt,line cap=round] (1,1) -- (3,0);
				\draw[lightgray,line width=\f*1.882pt,line cap=round] (3,0) -- (1,-1);
				\draw[lightgray,line width=\f*1.882pt,line cap=round] (-1,1) -- (-3,0);
				\draw[lightgray,line width=\f*1.882pt,line cap=round] (-2,0) -- (-3,0);
				\draw[lightgray,line width=\f*1.882pt,line cap=round] (-1,-1) -- (-3,0);
				\draw[darkgray,line width=\f*4.705pt,line cap=round] (1,1) -- (4,0);
				\draw[darkgray,line width=\f*4.705pt,line cap=round] (4,0) -- (1,-1);
				\draw[darkgray,line width=\f*4.705pt,line cap=round] (-1,1) -- (-4,0);
				\draw[darkgray,line width=\f*4.705pt,line cap=round] (-1,-1) -- (-4,0);
				\draw[lightgray,line width=\f*1.882pt,line cap=round] (3,0) -- (4,0);
				\draw[lightgray,line width=\f*1.882pt,line cap=round] (-3,0) -- (-4,0);
				\draw[lightgray,line width=\f*1.882pt,line cap=round] (4,0) -- (5,0);
				\draw[lightgray,line width=\f*1.882pt,line cap=round] (-4,0) -- (-5,0);
			\end{scope}
			\begin{scope}
				\filldraw[fill=white,draw=black,line width=\f*2pt] (0,0) circle[radius=\f*15pt];
				\node[anchor=north west, inner sep=1pt, font=\tiny] at (0,0) {$0$};
				\filldraw[fill=black, draw=black, line width=\f*2pt] (0,1) circle[radius=\f*15pt];
				\filldraw[fill=white,draw=black,line width=\f*2pt] (0,2) circle[radius=\f*15pt];
				\filldraw[fill=black, draw=black, line width=\f*2pt] (0,3) circle[radius=\f*15pt];
				\filldraw[fill=black, draw=black, line width=\f*2pt] (1,0) circle[radius=\f*15pt];
				\filldraw[fill=black, draw=black, line width=\f*2pt] (1,1) circle[radius=\f*15pt];
				\filldraw[fill=white,draw=black,line width=\f*2pt] (2,0) circle[radius=\f*15pt];
				\filldraw[fill=black, draw=black, line width=\f*2pt] (3,0) circle[radius=\f*15pt];
				\filldraw[fill=white,draw=black,line width=\f*2pt] (4,0) circle[radius=\f*15pt];
				\filldraw[fill=black, draw=black, line width=\f*2pt] (5,0) circle[radius=\f*15pt];
				\filldraw[fill=white,draw=black,line width=\f*2pt] (0,-1) circle[radius=\f*15pt];
				\filldraw[fill=white,draw=black,line width=\f*2pt] (0,-2) circle[radius=\f*15pt];
				\filldraw[fill=white,draw=black,line width=\f*2pt] (0,-3) circle[radius=\f*15pt];
				\filldraw[fill=white,draw=black,line width=\f*2pt] (1,-1) circle[radius=\f*15pt];
				\filldraw[fill=white,draw=black,line width=\f*2pt] (-1,0) circle[radius=\f*15pt];
				\filldraw[fill=black, draw=black, line width=\f*2pt] (-1,-1) circle[radius=\f*15pt];
				\filldraw[fill=white,draw=black,line width=\f*2pt] (-1,1) circle[radius=\f*15pt];
				\filldraw[fill=white,draw=black,line width=\f*2pt] (-2,0) circle[radius=\f*15pt];
				\filldraw[fill=white,draw=black,line width=\f*2pt] (-3,0) circle[radius=\f*15pt];
				\filldraw[fill=white,draw=black,line width=\f*2pt] (-4,0) circle[radius=\f*15pt];
				\filldraw[fill=white,draw=black,line width=\f*2pt] (-5,0) circle[radius=\f*15pt];
			\end{scope}
			\begin{scope}
				\draw[black,line width=\f*3.764pt,line cap=round,line join=round] (0.5,2) -- (0.5,1.5) -- (0,1.5) -- (-0.5,1) -- (0,0.5) -- (0.5,0.5) -- (0.5,0) -- (1,-0.5) -- (1.5,0) -- (1.5,0.5) -- (2.5,0) -- (2,-0.5) -- (3.5,0) -- (2.5,0.5) -- (3,0.5);
				\draw[black,line width=\f*3.764pt,line cap=round,line join=round] (-1.5,-1) -- (-2.5,-0.5) -- (-2,-0.5) -- (-1.5,-0.5) -- (-1,-0.5) -- (-0.5,-0.5) -- (-0.5,-1) -- (-0.5,-1.5) -- (-1,-1.5);
			\end{scope}
		\end{tikzpicture}
		\caption{$S$ is a two-edge double split.}
		\label{fig:double2}
	\end{subfigure}%
	\caption{The diamond zone $Z^\diamond$ for $\curve{\tri,\signof{\coll}}$ (top) and $\curve{\tri,\signof{\coll\setminus\{S\}}}$ (bottom) for different splits $S$. The respective curve is drawn in black, the split edges $S^\diamond$ in gray and the signs under $\signof{\coll}$ and $\signof{\coll\setminus\{S\}}$ of points in $Z^\diamond\cap A$ are indicated as black and white.}
	\label{fig:trivial_zones}
\end{figure}
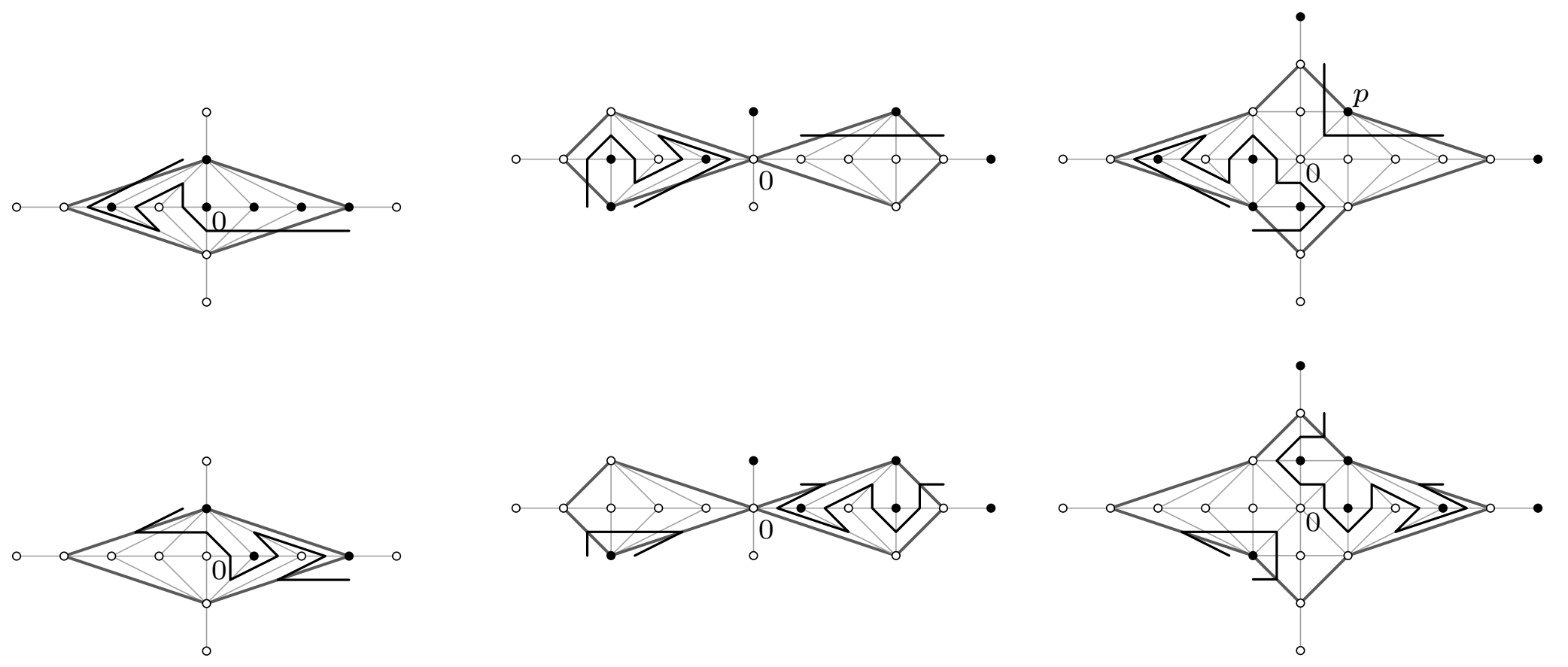

\begin{proposition}\label{prop:peel}
Let $\coll$ be a compatible collection of Harnack splits with at least one trivial member $S$.
Then $\curve{\coll}$ and $\curve{\coll\setminus\{S\}}$ have the same real scheme.
\end{proposition}
\begin{proof}
Let $Z$ be the trivial zone corresponding to $S$ and let $Z^\diamond$ be the \emph{diamond zone} of $Z$, i.e., the union of symmetric copies of $Z$ in the patchwork. Let $S^\diamond$ be the boundary of $Z^\diamond$. We investigate how the curves $\curve{\tri,\signof{\coll}}$ and $\curve{\tri,\signof{\coll\setminus\{S\}}}$ intersect $Z^\diamond$.
There are three possible cases based on $S$: 

\emph{Case (a):} $S$ is a simple split. In this case we may assume that $S$ has endpoints $(a,0)$ and $(0,b)$ with $a,b\neq 0$; otherwise apply a suitable map $g\in G$. If both $a$ and $b$ are greater than one, $Z$ contains the point $(1,1)$ and is hence not trivial. Thus we may assume $b=1$.
Note that $Z=Z^+$, so $\signof{\coll}$ and $\signof{\coll\setminus\{S\}}$ differ only on the interior of $Z$ and agree on $S\cap A$. 
The diamond zone $Z^\diamond$ has vertices $(a,0)$, $(0,1)$, $(-a,0)$, $(0,-1)$.
Observe that $\signof{\coll}(0,1)\neq\signof{\coll}(0,-1)$, hence $Z^\diamond$ has two opposite vertices with different signs.
Independently of the signs of $(a,0)$ and $(-a,0)$, this leaves only two edges whose endpoints have different signs (see Figure~\ref{fig:simple} for an example). Therefore 
$\curve{\tri,\signof{\coll}}$ crosses $S^\diamond$ in exactly these two edges. 
The same argument holds for $\curve{\tri,\signof{\coll\setminus\{S\}}}$. 
Since by Lemma~\ref{lem:ovals_contain_int_points} no oval can lie entirely within $Z^\diamond$, both curves intersect $Z^\diamond$ only in a single arc connecting the same points of $S^\diamond$.
Since furthermore $\curve{\tri,\signof{\coll}}$ and $\curve{\tri,\signof{\coll\setminus\{S\}}}$ coincide outside of $Z^\diamond$, the two curves are ambient isotopic.

\emph{Case (b):} $S$ is a double split both of whose endpoints lie on the same edge of $\dDelta$. In this case, we may assume that its points in $A$ are $(a,0), (c,1), (b,0)$ with $b=a+2k$ for some $k\in \mathbb{Z}_{>0}$, so they have the same parity. 
In this case $Z^\diamond$ consists of two blocks with vertices $(a,0)$, $(c,1)$, $(b,0)$, $(c,-1)$ and $(-a,0)$, $(-c,1)$, $(-b,0)$, $(-c,-1)$ respectively. Note that they may touch if $a=0$. Both blocks behave like $Z^\diamond$ in the previous case (see Figure~\ref{fig:double1}), and analogous arguments show that also the curves $\curve{\tri,\signof{\coll}}$ and $\curve{\tri,\signof{\coll\setminus\{S\}}}$ have the same isotopy type. 

\emph{Case (c):} $S$ is a double split connecting two edges of $\dDelta$. Here we may assume its points are $(a,0), (1,1), (0,b)$. In this case both $a$ and $b$ are even, since the endpoints must have the same parity. 
Note that $Z=Z^-$ only for $a=b=d$, in which case we may consider the opposite twist mask. So we may assume that $\signof{\coll}$ and $\signof{\coll\setminus\{S\}}$ differ only on the interior of $Z$. They again agree on $S\cap A$. 
The diamond zone $Z^\diamond$ has vertices $(a,0)$, $(1,1)$, $(0,b)$, $(-1,1)$, $(-a,0)$, $(-1,-1)$, $(0,-b)$, $(1,-1)$. 
Since $a$ and $b$ are both even, all points $(a,0),(0,b), (-a,0), (0,-b)$ have the same sign under $\signof{\coll}$, while $\signof{\coll}(1,1)=\signof{\coll}(-1,-1)\neq\signof{\coll}(1,-1)=\signof{\coll}(-1,1)$. So there are exactly two opposite points $p$ and $q$ with coordinates $\pm1$ whose sign differs from that of the other six points in $S^\diamond\cap A$.
Since $\signof{\coll}$ restricted to $Z$ is a Harnack sign function, one of these points, say $p$, only has neighbors with different signs within $Z^\diamond$; see Figure~\ref{fig:double2}. Therefore the arc intersecting $S^\diamond$ in one of the edges incident to $p$ cannot cross the coordinate axes within $Z^\diamond$. Consequently, the other arc must also cross $S^\diamond$ twice in the same quadrant. 
Moving from $\curve{\tri,\signof{\coll}}$ to $\curve{\tri,\signof{\coll\setminus\{S\}}}$ switches the roles of the two arcs, but the same argument holds. So in both cases the points where the curve crosses $S^\diamond$ are paired up in the same way. Since again no oval can lie entirely within $Z^\diamond$, the isotopy type is the same for $\curve{\tri,\signof{\coll}}$ and $\curve{\tri,\signof{\coll\setminus\{S\}}}$. 
\end{proof}

\begin{corollary}\label{cor:trivial}
Let $\cL\subseteq\splits$ be a list, and let $\cL_1$ be $\cL$ without its trivial members.
Then $\schemes{\cL_1}=\schemes{\cL}$.
\end{corollary}
\begin{proof}
  This follows by repeatedly removing trivial splits via Proposition~\ref{prop:peel}.
\end{proof}

Note that in Example~\ref{ex:degree2} all splits are trivial, so $\cL_1=\emptyset$ and only one real scheme can occur in degree~$2$.

\section{Results}\label{sec:results}

The backtracking walk of Section~\ref{sec:method} was implemented in \texttt{C++} and optimized for NVIDIA graphics cards.
For this, the isotopy algorithm of \cite{2026:GeiselmannEtAl:IsotopyComputation} was ported to CUDA in the library libisotopy, which is available at
\begin{center}
\libisotopyurl
\end{center}
The walk runs in a master and worker setup, in which the master process distributes subtrees of the search tree to the workers.
The code and the results, together with everything needed to reproduce them, are available at
\begin{center}
\splitenumurl
\end{center}

\thmdegreeeight*

\begin{proof}
By Corollary~\ref{cor:walk}, the result set of the walk contains all these schemes.
The run walked $873$ Harnack splits, \num{127764} compatible pairs and \num{105949232275} compatible collections up to the symmetry of Proposition~\ref{prop:symmetry}, the largest collection having $25$ splits.
It ran on 2026-10-26 on a single node with one AMD EPYC 7302P processor of $16$ cores and three NVIDIA GeForce RTX 2080 Ti cards, all $16$ cores and all three cards in use, and took \num{1367} seconds of wall time, about $23$ minutes, that is $6.1$ CPU core hours and $1.1$ GPU hours.
It reported no error and returned $38$ schemes.
Since the walk also encounters nonregular triangulations, we checked afterwards that each scheme is realized by a regular compatible collection, i.e., that the triangulation of some witness collection is regular, by solving a linear program \cite[\S5.2.1]{2010:DeLoeraRambauSantos:TriangulationsAlgorithms}.
Hence all $38$ schemes are realized by T-curves.
\end{proof}

\thmdegreenine*

\begin{proof}
By Corollary~\ref{cor:walk}, the result set of the walk contains all these schemes.
The run walked \num{1437} Harnack splits, \num{352227} compatible pairs and \num{57065860838381} compatible collections up to symmetry, the largest collection having $33$ splits.
It was run as a network of one master process to distribute the jobs, and $14$ worker nodes.
Four workers had one AMD EPYC 9254 processor and three NVIDIA RTX 6000 Ada Generation cards, nine had one AMD EPYC 7302 processor and three NVIDIA RTX A6000 cards, and one had one AMD EPYC 7302P processor and three NVIDIA Quadro RTX 6000 cards, so $224$ CPU cores and $42$ GPUs walked the tree.
The run took \num{15384} seconds of wall time, about $4.3$ hours, of which the workers were connected for \num{15244} seconds, that is $178$ GPU hours and $948$ allocated CPU core hours.
It returned $225$ schemes.
As in degree~$8$, each of them is realized by a regular compatible collection, hence by a T-curve.
\end{proof}

\section*{Acknowledgments}
This work was funded by the Deutsche Forschungsgemeinschaft under Germany's Excellence Strategy through the Berlin Mathematics Research Center MATH$^+$ (EXC-2046/1 and EXC-2046/2, project ID 390685689), through the Collaborative Research Centre TRR 195 Symbolic Tools in Mathematics and their Application, project ID 286237555, and under the National Research Data Infrastructure - NFDI 29/2 through the Mathematical Research Data Initiative MaRDI, project ID 460135501.

AI tools were used to help with coding and with editorial tasks, namely writing the abstract, correcting language and grammar, and maintaining overall consistency.
None of the proofs were done by AI.
The methods and all other creative input are the work of the human authors and were developed without the help of AI.

\printbibliography

\end{document}